\documentclass{article}
\usepackage{amssymb}
\usepackage{amsmath,amsfonts}

\begin{document}

\begin{center}
\bigskip {\Large Some remarks about Fibonacci elements in an arbitrary
algebra}%
\begin{equation*}
\end{equation*}

\bigskip Cristina FLAUT and Vitalii SHPAKIVSKYI

\begin{equation*}
\end{equation*}
\end{center}

\textbf{Abstract.  }{\small In this paper, we prove some relations between
Fibonacci elements in an arbitrary algebra. Moreover, we define imaginary
Fibonacci quaternions and imaginary Fibonacci octonions and we prove that
always three arbitrary imaginary Fibonacci quaternions are linear
independents and the mixed product of three arbitrary imaginary Fibonacci
octonions is zero.} 
\begin{equation*}
\end{equation*}



\begin{equation*}
\end{equation*}

\textbf{1. Introduction}%
\begin{equation*}
\end{equation*}

Fibonacci elements over some special algebras were intensively studied in
the last time in various papers, as for example: [Akk; ], [Fl, Sa; 15], [Fl,
Sh; 13(1)], [Fl, Sh; 13(2)], [Ha; ],[Ha1; ],[Ho; 61],[Ho; 63],[Ke; ]. All
these papers studied properties of Fibonacci quaternions, Fibonacci
octonions in Quaternion or Octonion algebras or in generalized Quaternion or
Octonion algebras or studied dual vectors or dual Fibonacci quaternions (
see [Gu;], [Nu; ]).

In this paper, we will prove that some of the obtained identities can be
obtained over an arbitrary algebras. We introduce the notions of imaginary
Fibonacci quaternions and imaginary Fibonacci octonions and we prove, using
the structure of the quaternion algebras \ and octonion algebras, that
always arbitrary three of such elements are linear dependents. For other
details, properties and applications regarding quaternion algebras \ and
octonion algebras, the reader is referred, for example, to [Sc; 54], [Sc;
66], [Fl, St; 09], [Sa, Fl, Ci; 09]. \bigskip

\textbf{2. Fibonacci elements in an arbitrary algebra}%
\begin{equation*}
\end{equation*}

\bigskip Let $A$ be a unitary algebra over $K$ ($K=\mathbb{R},\mathbb{C}$)
with a basis $\{e_{0}=1,e_{1},e_{2},...,e_{n}\}.$ Let $\{f_{n}\}_{n\in 
\mathbb{N}}$ be the Fibonacci sequence 
\begin{equation*}
f_{n}=f_{n-1}+f_{n-2},n\geq 2,f_{0}=0,f_{1}=1.
\end{equation*}%
In algebra $A,$ we define the Fibonacci element as follows:%
\begin{equation*}
F_{m}=\overset{n}{\underset{k=0}{\sum }}f_{m+k}e_{k}.
\end{equation*}

\textbf{Proposition 2.1.} \textit{With the above notations, the following
relations hold:}

\textit{1)} $F_{m+2}=F_{m+1}+F_{m};$

\textit{2)} $\overset{p}{\underset{i=0}{\sum }}F_{i}=F_{p+2}-F_{i}.$

\textbf{Proof.} 1) $F_{m+1}+F_{m}=\overset{n}{\underset{k=0}{\sum }}%
f_{m+k+1}e_{k}+\overset{n}{\underset{k=0}{\sum }}f_{m+k}e_{k}=\overset{n}{%
\underset{k=0}{\sum }}(f_{m+k+1}+f_{m+k})e_{k}=\overset{n}{\underset{k=0}{%
\sum }}$ $f_{m+k+2}e_{k}=F_{m+2}.$

2) $\overset{p}{\underset{i=0}{\sum }}F_{i}=F_{1}+F_{2}+...+F_{p}=$\newline
$=\overset{n}{\underset{k=0}{\sum }}f_{k+1}e_{k}+\overset{n}{\underset{k=0}{%
\sum }}f_{k+2}e_{k}+...+\overset{n}{\underset{k=0}{\sum }}f_{k+p}e_{k}=$%
\newline
$=e_{0}\left( f_{1}+...+f_{p}\right) +e_{1}\left( f_{2}+...+f_{p+1}\right) +$%
\newline
$+e_{2}\left( f_{3}+...+f_{p+2}\right) +...+e_{n}\left(
f_{k+n}+...+f_{p+n}\right) =$\newline
$=$ $e_{0}\left( f_{p+2}-1\right) +e_{1}\left( f_{p+3}-1-f_{1}\right)
+e_{2}\left( f_{p+4}-1-f_{1}-f_{2}\right) +$\newline
$+e_{3}\left( f_{p+5}-1-f_{1}-f_{2}-f_{3}\right) +...+e_{n}\left(
f_{p+n+2}-1-f_{1}-f_{2}-...-f_{n}\right) =$\newline
$=F_{p+2}-F_{2}.$

We used the identity $\overset{p}{\underset{i=1}{\sum }}f_{i}=f_{p+2}-1$
(for usual Fibonacci numbers) and $1+f_{1}+f_{2}+...+f_{n}=f_{n+2}.\Box
\medskip $

\textbf{Remark 2.2}. The equalities 1, 2 from the above proposition
generalize the corresponding formulae from [Ke; ] [Ha; ] [Nu; ] [Ha1;
].\bigskip

\textbf{Proposition 2.3.} \textit{We have the following formula (Binet's
fomula):} 
\begin{equation*}
F_{m}=\frac{\alpha ^{\ast }\alpha ^{m}-\beta ^{\ast }\beta ^{m}}{\alpha
-\beta },
\end{equation*}%
\textit{where} $\alpha =\frac{1+\sqrt{5}}{2},\beta =\frac{1-\sqrt{5}}{2}%
,\alpha ^{\ast }=\underset{k=0}{\overset{n}{\sum }}\alpha ^{k}e_{k},~\beta
^{\ast }=\underset{k=0}{\overset{n}{\sum }}\beta ^{k}e_{k}.$

\textbf{Proof.} \ Using the formula for the real quaternions, $f_{m}=\frac{%
\alpha ^{m}-\beta ^{m}}{\alpha -\beta },$ we obtain\newline
$F_{m}=\overset{n}{\underset{k=0}{\sum }}f_{m+k}e_{k}=\frac{\alpha
^{m}-\beta ^{m}}{\alpha -\beta }e_{0}+\frac{\alpha ^{m+1}-\beta ^{m+1}}{%
\alpha -\beta }e_{1}+\frac{\alpha ^{m+2}-\beta ^{m+2}}{\alpha -\beta }%
e_{2}+...+$\newline
$+\frac{\alpha ^{m+n}-\beta ^{m+n}}{\alpha -\beta }e_{n}=\frac{a^{m}}{\alpha
-\beta }\left( e_{0}+\alpha e_{1}+\alpha ^{2}e_{2}+...+\alpha
^{n}e_{n}\right) +$\newline
$+\frac{\beta ^{m}}{\alpha -\beta }\left( e_{0}+\beta e_{1}+\beta
^{2}e_{2}+...+\beta ^{n}e_{n}\right) =\frac{\alpha ^{\ast }\alpha ^{m}-\beta
^{\ast }\beta ^{m}}{\alpha -\beta }.\medskip \Box \medskip $

\textbf{Remark 2.4.} The above result generalizes the Binet formulae from
the papers [Gu;] [Akk; ] [Ke; ] [Ha; ] [Nu; ] [Ha1; ].\medskip

\textbf{Theorem 2.5.} \textit{The generating function for the Fibonacci
number over an algebra is of the form} 
\begin{equation*}
G\left( t\right) =\frac{F_{0}+\left( F_{1}-F_{0}\right) t}{1-t-t^{2}}.
\end{equation*}

\textbf{Proof.} We consider the generating function of the form%
\begin{equation*}
G\left( t\right) =\overset{\infty }{\underset{m=0}{\sum }}F_{m}t^{m}.
\end{equation*}%
We consider the product\newline
$G\left( t\right) \left( 1-t-t^{2}\right) =\overset{\infty }{\underset{m=0}{%
\sum }}F_{m}t^{m}=$ $\overset{\infty }{\underset{m=0}{\sum }}F_{m}t^{m}-%
\overset{\infty }{\underset{m=0}{\sum }}F_{m}t^{m+1}-\overset{\infty }{%
\underset{m=0}{\sum }}F_{m}t^{m+2}=$\newline
$=F_{0}+F_{1}t+F_{2}t^{2}+F_{3}t^{3}+...-F_{0}t-F_{1}t^{2}-F_{2}t^{3}-...-$%
\newline
$-F_{0}t^{2}-F_{1}t^{3}-F_{2}t^{4}-...=F_{0}+\left( F_{1}-F_{0}\right)
t.\Box \medskip $

\textbf{Remark 2.6.} The above Theorem generalizes results from the papers
[Gu;], [Akk; ], \ [Ke; ], [Ha; ],[Nu; ].

\begin{equation*}
\end{equation*}

\textbf{The Cassini identity}%
\begin{equation*}
\end{equation*}

First, we obtain the following identity.

\textbf{Proposition 2.7.}

\begin{equation}
F_{-m}=\left( -1\right) ^{m+1}f_{m}F_{1}+\left( -1\right) ^{m}f_{m+1}F_{0}. 
\tag{2.1.}
\end{equation}

\textbf{Proof.} We use induction. \ For $m=1,$ we obtain $%
F_{-1}=f_{1}F_{1}-f_{2}F_{0},$ which is true. Now, we assume that it is true
for an arbitrary integer $k$%
\begin{equation*}
F_{-k}=\left( -1\right) ^{k+1}f_{k}F_{1}+\left( -1\right) ^{k}f_{k+1}F_{0}
\end{equation*}%
For $k+1,$ we obtain\newline
$F_{-(k+1)}=\left( -1\right) ^{k+2}f_{k+1}F_{1}+\left( -1\right)
^{k+1}f_{k+2}F_{0}=$\newline
$=\left( -1\right) ^{k}f_{k}F_{1}+\left( -1\right) ^{k}f_{k-1}F_{1}+\left(
-1\right) ^{k-1}f_{k+1}F_{0}+$\newline
$+\left( -1\right) ^{k-1}f_{k}F_{0}=F_{-\left( n-1\right) }-F_{-n}.$
Therefore, this statement is true.$\Box \medskip $

\textbf{Theorem 2.8.} (Cassini's identity) \textit{With the above notations,
we have the following formula}

\begin{equation*}
F_{m-1}F_{m+1}-F_{m}^{2}=\left( -1\right) ^{m}(F_{-1}F_{1}-F_{0}^{2}).
\end{equation*}

\textbf{Proof.}

We consider\newline
$%
F_{m-1}=f_{m-1}e_{0}+f_{m}e_{1}+f_{m+1}e_{2}+f_{m+2}e_{3}+...+f_{m+n-1}e_{n}, 
$\newline
$%
F_{m+1}=f_{m+1}e_{0}+f_{m+2}e_{1}+f_{m+3}e_{2}+f_{m+4}e_{3}+...+f_{m+n+1}e_{n}, 
$\newline
$F_{m}=f_{m}e_{0}+f_{m+1}e_{1}+f_{m+2}e_{2}+f_{m+2}e_{3}+...+f_{m+n}e_{n}.$

We compute\newline
$F_{m-1}F_{m+1}=$\newline
$=\left[
f_{m-1}f_{m+1}e_{0}^{2}+f_{m-1}f_{m+2}e_{0}e_{1}+f_{m-1}f_{m+3}e_{0}e_{2}+f_{m-1}f_{m+4}e_{0}e_{3}...+f_{m-1}f_{m+n+1}e_{0}e_{n}%
\right] +$\newline
$%
+[f_{m}f_{m+1}e_{1}e_{0}+f_{m}f_{m+2}e_{1}^{2}+f_{m}f_{m+3}e_{1}e_{2}+f_{m}f_{m+4}e_{1}e_{3}...+f_{m}f_{m+n+1}e_{1}e_{n}]+ 
$\newline
$%
+[f_{m+1}^{2}e_{2}e_{0}+f_{m+1}f_{m+2}e_{2}e_{1}+f_{m+1}f_{m+3}e_{2}^{2}+f_{m+1}f_{m+4}e_{1}e_{3}...+f_{m+1}f_{m+n+1}e_{2}e_{n}]+ 
$\newline
$%
+[f_{m+2}f_{m+1}e_{3}e_{0}+f_{m+2}^{2}e_{3}e_{1}+f_{m+2}f_{m+3}e_{3}e_{2}+f_{m+2}f_{m+4}e_{3}^{2}...+f_{m+2}f_{m+n+1}e_{3}e_{n}]+...+ 
$\newline
$%
+[f_{m+n-1}f_{m+1}e_{n}e_{0}+f_{m+n-1}f_{m+2}e_{n}e_{1}+f_{m+n-1}f_{m+3}e_{n}e_{2}+f_{m+n-1}f_{m+4}e_{n}e_{3}...+f_{m+n-1}f_{m+n+1}e_{n}^{2}]. 
$

Now, we compute \newline
$F_{m}^{2}=\left[
f_{m}^{2}e_{0}^{2}+f_{m}f_{m+1}e_{0}e_{1}+f_{m}f_{m+2}e_{0}e_{2}+f_{m}f_{m+3}e_{0}e_{3}+...+f_{m}f_{m+n}e_{0}e_{n}%
\right] +$\newline
$+$ $\left[
f_{m+1}f_{m}e_{1}e_{0}+f_{m+1}^{2}e_{1}^{2}+f_{m+1}f_{m+2}e_{1}e_{2}+f_{m+1}f_{m+3}e_{1}e_{3}+...+f_{m+1}f_{m+n}e_{1}e_{n}%
\right] +$\newline
$+\left[
f_{m+2}f_{m}e_{2}e_{0}+f_{m+2}f_{m+1}e_{2}e_{1}+f_{m+2}^{2}e_{2}^{2}+f_{m+2}f_{m+3}e_{2}e_{3}+...+f_{m+2}f_{m+n}e_{2}e_{n}%
\right] +$\newline
$+\left[
f_{m+2}f_{m}e_{2}e_{0}+f_{m+2}f_{m+1}e_{2}e_{1}+f_{m+2}^{2}e_{2}^{2}+f_{m+2}f_{m+3}e_{2}e_{3}+...+f_{m+2}f_{m+n}e_{2}e_{n}%
\right] +$\newline
$+\left[
f_{m+3}f_{m}e_{3}e_{0}+f_{m+3}f_{m+1}e_{3}e_{1}+f_{m+3}f_{m+2}e_{3}e_{2}+f_{m+3}^{2}e_{3}^{2}+...+f_{m+3}f_{m+n}e_{3}e_{n}%
\right] +...+$\newline
$+\left[
f_{m+n}f_{m}e_{n}e_{0}+f_{m+n}f_{m+1}e_{n}e_{1}+f_{m+n}f_{m+2}e_{n}e_{2}+f_{m+n}f_{m+3}e_{n}e_{3}+...+f_{m+n}^{2}e_{n}^{2}%
\right] .$

Consider the difference\newline
$F_{m-1}F_{m+1}-F_{m}^{2}=$\newline
$=e_{0}\left[ e_{0}\left( f_{m-1}f_{m+1}-f_{m}^{2}\right) +e_{1}\left(
f_{m-1}f_{m+2}-f_{m}f_{m+1}\right) +...+e_{n}\left(
f_{m-1}f_{m+n+1}-f_{m}f_{m+n}\right) \right] +$\newline
$+e_{1}\left[ e_{0}\left( f_{m}f_{m+1}-f_{m+1}f_{m}\right) +e_{1}\left(
f_{m}f_{m+2}-f_{m+1}^{2}\right) +...+e_{n}\left(
f_{m}f_{m+n+1}-f_{m+1}f_{m+n}\right) \right] +$\newline
$+e_{2}\left[ e_{0}\left( f_{m+1}^{2}-f_{m+2}f_{m}\right) +e_{1}\left(
f_{m+1}f_{m+2}-f_{m+2}f_{m+1}\right) +...+e_{n}\left(
f_{m+1}f_{m+n+1}-f_{m+2}f_{m+n}\right) \right] +$\newline
$+e_{3}\left[ e_{0}\left( f_{m+2}f_{m+1}-f_{m+3}f_{m}\right) +e_{1}\left(
f_{m+2}^{2}-f_{m+3}f_{m+1}\right) +...+e_{n}\left(
f_{m+2}f_{m+n+1}-f_{m+3}f_{m+n}\right) \right] +...+$\newline
$+e_{n}\left[ e_{0}\left( f_{m+n-1}f_{m+1}-f_{m+n}f_{m}\right) +e_{1}\left(
f_{m+n-1}f_{m+2}-f_{m+n}f_{m+1}\right) +...+e_{n}\left(
f_{m+n-1}f_{m+n+1}-f_{m+n}^{2}\right) \right] .$

Using the formula $f_{i}f_{j}-f_{i+k}f_{j-k}=\left( -1\right)
^{j-k}f_{i+k-j}f_{k}$ (see Koshy, p. 87, formula 2) and the identities $%
f_{1}=1,f_{-m}=\left( -1\right) ^{m+1}f_{m}$ (see Koshy, p. 84), we obtain%
\newline
$F_{m-1}F_{m+1}-F_{m}^{2}=e_{0}\left( -1\right) ^{m+1}\left[
e_{0}f_{1}+e_{1}f_{2}+e_{2}f_{3}+...+e_{n}f_{n+1}\right] +$\newline
$+e_{1}\left( -1\right) ^{m+1}\left[
e_{0}f_{0}+e_{1}f_{1}+e_{2}f_{2}+...+e_{n}f_{n}\right] +$\newline
$+e_{2}\left( -1\right) ^{m}\left[
e_{0}f_{-1}+e_{1}f_{0}+e_{2}f_{1}+...+e_{n}f_{n-1}\right] +$\newline
$+e_{3}\left( -1\right) ^{m}\left[
e_{0}f_{-2}+e_{1}f_{-1}+e_{2}f_{0}+...+e_{n}f_{n-2}\right] +...+$\newline
$+\left( -1\right) ^{m\ast n}e_{n}\left[
e_{0}f_{-n+1}+e_{1}f_{-n+2}+e_{2}f_{-n+3}+...+e_{n}f_{1}\right] =$\newline
$=\left( -1\right) ^{m}\left(
e_{0}F_{1}-e_{1}F_{0}+e_{2}F_{-1}-e_{3}F_{-2}+...+\left( -1\right)
^{n}e_{n}F_{-n+1}\right) .$

Using Proposition 2.7, we have\newline
$F_{m-1}F_{m+1}-F_{m}^{2}=\left( -1\right)
^{m}[e_{0}F_{1}-e_{1}F_{0}+e_{2}\left( F_{1}-F_{0}\right) -$\newline
$-e_{3}\left( 2F_{0}-F_{1}\right) +e_{4}\left( 2F_{1}-3F_{0}\right)
-e_{5}\left( -3F_{1}+5F_{0}\right) +...+$\newline
$+e_{n}\left( -1\right) ^{n}\left( \left( -1\right) ^{n}f_{n-1}F_{1}+\left(
-1\right) ^{n-1}f_{n}F_{0}\right) ]=$\newline
$=\left( -1\right)
^{m}[(e_{0}f_{-1}+e_{1}f_{0}+e_{2}f_{1}+...+e_{n}f_{n-1})F_{1}-$\newline
$-\left( f_{0}e_{0}+f_{1}e_{1}+f_{2}e_{2}+...+f_{n}e_{n}\right) F_{0}]=$%
\newline
$=\left( -1\right) ^{m}\left[ F_{-1}F_{1}-F_{0}^{2}\right] .$ The theorem is
now proved.\medskip

\textbf{Remark 2.9. }

i) Similarly, we can prove an analogue of Cassini's formula:%
\begin{equation*}
F_{m+1}F_{m-1}-F_{m}^{2}=\left( -1\right) ^{m}\left[ F_{1}F_{-1}-F_{0}^{2}%
\right] .
\end{equation*}

ii) Theorem 2.8 generalizes Cassini's formula for all real algebras.

iii) If the algebra $A$ is algebra of the real numbers $\mathbb{R},$ in this
case, we have $F_{m}=f_{m}.$ From the above theorem, it \ results that

\begin{equation*}
f_{m+1}f_{m-1}-f_{m}^{2}=\left( -1\right) ^{m}\left[ f_{1}f_{-1}-f_{0}^{2}%
\right] =\left( -1\right) ^{m},
\end{equation*}%
which it is the classical Cassini's identity.%
\begin{equation*}
\end{equation*}

\bigskip \textbf{3. Imaginary Fibonacci quaternions and imaginary Fibonacci
octonions }%
\begin{equation*}
\end{equation*}

In the following, we will consider a field $K$ with $charK\neq 2,3,$ $V$ a
finite dimensional vector space and $A$ a finite dimensional unitary algebra
over a field $\ K$, associative or nonassociative.

Let $\mathbb{H}\left( \alpha ,\beta \right) $ be the generalized real\
quaternion algebra, the algebra of the elements of the form $a=a_{1}\cdot
1+a_{2}\mathbf{i}+a_{3}\mathit{j}+a_{4}\mathbf{k},$ where $a_{i}\in \mathbb{R%
},\mathbf{i}^{2}=-\alpha ,\mathbf{j}^{2}=-\beta ,$ $\mathbf{k}=\mathbf{ij}=-%
\mathbf{ji}.$ We denote by $\mathbf{t}\left( a\right) $ and $\mathbf{n}%
\left( a\right) $ the trace and the norm of a real quaternion $a.$ The norm
of a generalized quaternion has the following expression $\mathbf{n}\left(
a\right) =a_{1}^{2}+\alpha a_{2}^{2}+\beta a_{3}^{2}+\alpha \beta a_{4}^{2}$
and $\mathbf{t}\left( a\right) =2a_{1}.$ It is known that for $a\in $ $%
\mathbb{H}\left( \alpha ,\beta \right) ,$ we have $a^{2}-\mathbf{t}\left(
a\right) a+\mathbf{n}\left( a\right) =0.$ The quaternion algebra $\mathbb{H}%
\left( \alpha ,\beta \right) $ is a \textit{division algebra} if for all $%
a\in \mathbb{H}\left( \alpha ,\beta \right) ,$ $a\neq 0$ we have $\mathbf{n}%
\left( a\right) \neq 0,$ otherwise $\mathbb{H}\left( \alpha ,\beta \right) $
is called a \textit{split algebra}.

Let $\mathbb{O}(\alpha ,\beta ,\gamma )$ be a generalized octonion algebra
over $\mathbb{R},$ with basis $\{1,e_{1},...,e_{7}\},$ the algebra of the
elements of the form $\
a=a_{0}+a_{1}e_{1}+a_{2}e_{2}+a_{3}e_{3}+a_{4}e_{4}+a_{5}e_{5}+a_{6}e_{6}+a_{7}e_{7}\, 
$and the multiplication given in the following table:\medskip \medskip

\begin{center}
{\footnotesize $%
\begin{tabular}{c||c|c|c|c|c|c|c|c|}
$\cdot $ & $1$ & $\,\,\,e_{1}$ & $\,\,\,\,\,e_{2}$ & $\,\,\,\,e_{3}$ & $%
\,\,\,\,e_{4}$ & $\,\,\,\,\,\,e_{5}$ & $\,\,\,\,\,\,e_{6}$ & $%
\,\,\,\,\,\,\,e_{7}$ \\ \hline\hline
$\,1$ & $1$ & $\,\,\,e_{1}$ & $\,\,\,\,e_{2}$ & $\,\,\,\,e_{3}$ & $%
\,\,\,\,e_{4}$ & $\,\,\,\,\,\,e_{5}$ & $\,\,\,\,\,e_{6}$ & $%
\,\,\,\,\,\,\,e_{7}$ \\ \hline
$\,e_{1}$ & $\,\,e_{1}$ & $-\alpha $ & $\,\,\,\,e_{3}$ & $-\alpha e_{2}$ & $%
\,\,\,\,e_{5}$ & $-\alpha e_{4}$ & $-\,\,e_{7}$ & $\,\,\,\alpha e_{6}$ \\ 
\hline
$\,e_{2}$ & $\,e_{2}$ & $-e_{3}$ & $-\,\beta $ & $\,\,\beta e_{1}$ & $%
\,\,\,\,e_{6}$ & $\,\,\,\,\,e_{7}$ & $-\beta e_{4}$ & $-\beta e_{5}$ \\ 
\hline
$e_{3}$ & $e_{3}$ & $\alpha e_{2}$ & $-\beta e_{1}$ & $-\alpha \beta $ & $%
\,\,\,\,e_{7}$ & $-\alpha e_{6}$ & $\,\,\,\beta e_{5}$ & $-\alpha \beta
e_{4} $ \\ \hline
$e_{4}$ & $e_{4}$ & $-e_{5}$ & $-\,e_{6}$ & $-\,\,e_{7}$ & $-\,\gamma $ & $%
\,\,\,\gamma e_{1}$ & $\,\,\gamma e_{2}$ & $\,\,\,\,\,\gamma e_{3}$ \\ \hline
$\,e_{5}$ & $\,e_{5}$ & $\alpha e_{4}$ & $-\,e_{7}$ & $\,\alpha e_{6}$ & $%
-\gamma e_{1}$ & $-\,\alpha \gamma $ & $-\gamma e_{3}$ & $\,\alpha \gamma
e_{2}$ \\ \hline
$\,\,e_{6}$ & $\,\,e_{6}$ & $\,\,\,\,e_{7}$ & $\,\,\beta e_{4}$ & $-\,\beta
e_{5}$ & $-\gamma e_{2}$ & $\,\,\,\gamma e_{3}$ & $-\beta \gamma $ & $-\beta
\gamma e_{1}$ \\ \hline
$\,\,e_{7}$ & $\,\,e_{7}$ & $-\alpha e_{6}$ & $\,\beta e_{5}$ & $\alpha
\beta e_{4}$ & $-\gamma e_{3}$ & $-\alpha \gamma e_{2}$ & $\beta \gamma
e_{1} $ & $-\alpha \beta \gamma $ \\ \hline
\end{tabular}%
\ \medskip $ }\medskip

Table 1
\end{center}

The algebra $\mathbb{O}(\alpha ,\beta ,\gamma )$ is non-commutative and
non-associative.

If $a\in \mathbb{O}(\alpha ,\beta ,\gamma ),$ $%
a=a_{0}+a_{1}e_{1}+a_{2}e_{2}+a_{3}e_{3}+a_{4}e_{4}+a_{5}e_{5}+a_{6}e_{6}+a_{7}e_{7} 
$ then $\bar{a}%
=a_{0}-a_{1}e_{1}-a_{2}e_{2}-a_{3}e_{3}-a_{4}e_{4}-a_{5}e_{5}-a_{6}e_{6}-a_{7}e_{7} 
$ is called the \textit{conjugate} of the element $a.$ The scalars $\mathbf{t%
}\left( a\right) =a+\overline{a}\in \mathbb{R}$ and 
\begin{equation}
\,\mathbf{n}\left( a\right) =a\overline{a}=a_{0}^{2}+\alpha a_{1}^{2}+\beta
a_{2}^{2}+\alpha \beta a_{3}^{2}+\gamma a_{4}^{2}+\alpha \gamma
a_{5}^{2}+\beta \gamma a_{6}^{2}+\alpha \beta \gamma a_{7}^{2}\in \mathbb{R},
\tag{3.1.}
\end{equation}%
are called the \textit{trace}, respectively, the \textit{norm} of the
element $a\in $ $A.$ \thinspace It\thinspace \thinspace \thinspace
follows\thinspace \thinspace \thinspace that$\,$\newline
\thinspace \thinspace $a^{2}-\mathbf{t}\left( a\right) a+\mathbf{n}\left(
a\right) =0,\forall a\in A.$The octonion algebra $\mathbb{O}\left( \alpha
,\beta ,\gamma \right) $ is a \textit{division algebra} if for all $a\in 
\mathbb{O}\left( \alpha ,\beta ,\gamma \right) ,$ $a\neq 0$ we have $\mathbf{%
n}\left( a\right) \neq 0,$ otherwise $\mathbb{O}\left( \alpha ,\beta ,\gamma
\right) $ is called a \textit{split algebra}.\medskip

Let \ $V$ be a real vector space of dimension $n$ and $<,>$ be the inner
product. The\textit{\ cross product} on $V$ is a continuos map%
\begin{equation*}
X:V^{s}\rightarrow V,s\in \{1,2,...,n\}
\end{equation*}%
with the following properties:

1) $<X\left( x_{1},...x_{s}\right) ,x_{i}>=0,i\in \{1,2,...,s\};$

2) $<X\left( x_{1},...x_{s}\right) ,X\left( x_{1},...x_{s}\right) >=\det
\left( <x_{i},x_{j}>\right) .($see [Br; ]$)$

In [Ro; 96], was proved that if $d=\dim _{\mathbb{R}}V,$ therefore $d\in
\{0,1,3,7\}.$(see [Ro; 96], Proposition 3)

The values $0,1,3$ and $7$ for dimensions are obtained from Hurwitz's
theorem, since the real Hurwitz division algebras $\mathcal{H}$ exist only
for dimensions $1,2,4$ and $8$ dimensions. In this situations, the cross
product is obtained from the product of the normed division algebra,
restricting it to imaginary subspace of the algebra $\mathcal{H},$ which can
be of \ dimension $0,1,3$ or $7$.(see [Ja; 74]) It is known that the real
Hurwitz division algebras are only: the real numbers, the complex numbers,
the quaternions and the octonions.

In $\mathbb{R}^{3}$ with the canonical basis $\{i_{1},i_{2},i_{3}\},$ the
cross product of two linearly independent vectors $%
x=x_{1}i_{1}+x_{2}i_{2}+x_{3}i_{3}$ and $y=y_{1}i_{1}+y_{2}i_{2}+y_{3}i_{3}$
is a vector, denoted by $x\times y$ and can be expressed computing the
following formal determinant 
\begin{equation}
x\times y=\left\vert 
\begin{array}{ccc}
i_{1} & i_{2} & i_{3} \\ 
x_{1} & x_{2} & x_{3} \\ 
y_{1} & y_{2} & y_{3}%
\end{array}%
\right\vert .  \tag{3.2.}
\end{equation}

The cross product can also be described using the quaternions and the basis $%
\{i_{1},i_{2},i_{3}\}$ as a standard basis for $\mathbb{R}^{3}.$ If a vector 
$x\in \mathbb{R}^{3}$ has the form $x=x_{1}i_{1}+x_{2}i_{2}+x_{3}i_{3}$ and
is represented as the quaternion $x=x_{1}\mathbf{i}+x_{2}\mathbf{j}+x_{3}%
\mathbf{k}$, therefore the cross product of two vectors has the form $%
x\times y=xy+<x,y>,$ where $<x,y>=x_{1}y_{1}+x_{2}y_{2}+x_{3}y_{3}$ is the
inner product.

A cross product for 7-dimensional vectors can be obtained in the same way by
using the octonions instead of the quaternions. If $x=\underset{i=0}{\overset%
{7}{\sum }}x_{i}e_{i}$ and $y=\underset{i=0}{\overset{7}{\sum }}y_{i}e_{i}$
are two imaginary octonions, therefore%
\begin{eqnarray*}
x\times y
&=&(x_{2}y_{4}-x_{4}y_{2}+x_{3}y_{7}-x_{7}y_{3}+x_{5}y_{6}-x_{6}y_{5})%
\,e_{1}+ \\
&&+(x_{3}y_{5}-x_{5}y_{3}+x_{4}y_{1}-x_{1}y_{4}+x_{6}y_{7}-x_{7}y_{6})%
\,e_{2}+ \\
&&+(x_{4}y_{6}-x_{6}y_{4}+x_{5}y_{2}-x_{2}y_{5}+x_{7}y_{1}-x_{1}y_{7})%
\,e_{3}+ \\
&&+(x_{5}y_{7}-x_{7}y_{5}+x_{6}y_{3}-x_{3}y_{6}+x_{1}y_{2}-x_{2}y_{1})%
\,e_{4}+ \\
&&+(x_{6}y_{1}-x_{1}y_{6}+x_{7}y_{4}-x_{4}y_{7}+x_{2}y_{3}-x_{3}y_{2})%
\,e_{5}+ \\
&&+(x_{7}y_{2}-x_{2}y_{7}+x_{1}y_{5}-x_{5}y_{1}+x_{3}y_{4}-x_{4}y_{3})%
\,e_{6}+ \\
&&+(x_{1}y_{3}-x_{3}y_{1}+x_{2}y_{6}-x_{6}y_{2}+x_{4}y_{5}-x_{5}y_{4})%
\,e_{7},
\end{eqnarray*}%
\begin{equation}
\tag{3.3.}
\end{equation}%
see [Si; 02].

Let $\mathbb{H}$ be the real division quaternion algebra and $\mathbb{H}%
_{0}=\{x\in \mathbb{H}$ $\mid $ $\mathbf{t}\left( x\right) =0\}.$ An element 
$F_{n}\in \mathbb{H}_{0}$ is called an\textit{\ imaginary Fibonacci
quaternion element} if it is of the form $F_{n}=f_{n}\mathbf{i}+f_{n+1}%
\mathbf{j}+f_{n+2}\mathbf{k,}$ where $\left( f_{n}\right) _{n\in \mathbb{N}}$
is the Fibonacci numbers sequence$.$ Let $F_{k},F_{m},F_{n}$ be three
imaginary Fibonacci quaternions. Therefore, we have the following result.

In the proof of the following results, we will use some relations between
Fibonacci numbers, namely:

\textit{D'Ocagne's identity}%
\begin{equation}
f_{m}f_{n+1}-f_{n}f_{m+1}=\left( -1\right) ^{n}f_{m-n}  \tag{3.4.}
\end{equation}%
see relation (33) from [Wo], and

\textit{Johnson's identity}%
\begin{equation}
f_{a}f_{b}-f_{c}f_{d}=\left( -1\right) ^{r}\left(
f_{a-r}f_{b-r}-f_{c-r}f_{d-r}\right) ,  \tag{3.5.}
\end{equation}%
for arbitrary integers $a,b,c,d,$ and $r$ with $a+b=c+d,$ see relation (36)
from [Wo].

\textbf{Proposition 3.1.}  \textit{With the above notations, for three
arbitrary Fibonacci imaginary quaternions, } \textit{we have}%
\begin{equation*}
<F_{k}\times F_{m},F_{n}>=0.
\end{equation*}%
\textit{Therefore, the vectors} $F_{k},F_{m},F_{n}$ \textit{are linear
dependents}.\medskip 

The above result is similar with the result for dual Fibonacci vectors
obtained in [Gu;], Theorem 11.

Let $\mathbb{O}$ be the real division octonion algebra and $\mathbb{O}%
_{0}=\{x\in \mathbb{H}$ $\mid $ $\mathbf{t}\left( x\right) =0\}.$ An element 
$F_{n}\in \mathbb{O}_{0}$ is called an \textit{imaginary Fibonacci octonion
element} if it is of the form $%
F_{n}=f_{n}e_{1}+f_{n+1}e_{1}+f_{n+2}e_{1}+f_{n+3}e_{1}+f_{n+4}e_{1}+f_{n+5}e_{1}+f_{n+6}e_{1}%
\mathbf{,}$ where $\left( f_{n}\right) _{n\in \mathbb{N}}$ is the Fibonacci
numbers sequence$.$ Let $F_{k},F_{m},F_{n}$ be three imaginary Fibonacci
octonions.\medskip\ 

\textbf{Proposition 3.2.} \textit{With the above notations, for three
arbitrary Fibonacci imaginary octonions, we have}%
\begin{equation*}
<F_{k}\times F_{m},F_{n}>=0.
\end{equation*}%
\qquad 

\textbf{Proof.} Using formulae ($3.3$), $\left( 3.4\right) $ and $\left(
3.5\right) $, we will compute \medskip $F_{k}\times F_{m}.$ \newline
The coefficient of \ $e_{1}$ is\newline
$f_{m+2}f_{k+4}-f_{k+2}f_{m+4}+f_{m+3}f_{k+7}-f_{k+3}f_{m+7}+f_{m+5}f_{k+6}-$
$f_{k+5}f_{m+6}=$\newline
$=f_{m}f_{k+2}-f_{k}f_{m+2}-f_{m}f_{k+4}+f_{k}f_{m+4}-f_{m}f_{k+1}+$ $%
f_{k}f_{m+1}=$\newline
$=f_{m}\left( f_{k+2}-f_{k+4}-f_{k+1}\right) +f_{k}\left(
-f_{m+2}+f_{m+4}+f_{m+1}\right) =$\newline
$=f_{m}\left( f_{k}-f_{k+4}\right) +f_{k}\left( f_{m+4}-f_{m}\right) =$%
\newline
$=-f_{m}\left( 3f_{k+1}+f_{k}\right) +f_{k}\left( 3f_{m+1}+f_{m}\right) =$%
\newline
$=-3\left( f_{m}f_{k+1}-f_{k}f_{m+1}\right) =-3\left( -1\right) ^{k}f_{m-k}.$%
\newline
The coefficient of $e_{2}$ is\newline
$f_{m+3}f_{k+5}-f_{k+3}f_{m+5}+f_{m+4}f_{k+1}-f_{k+4}f_{m+1}+f_{m+6}f_{k+7}-$
$f_{k+6}f_{m+7}=$\newline
$=$ $-f_{m}f_{k+2}+f_{k}f_{m+2}-f_{m+3}f_{k}+f_{k+3}f_{m}+f_{m}f_{k+1}-$ $%
f_{k}f_{m+1}=$\newline
$=f_{m}\left( -f_{k+2}+f_{k+3}+f_{k+1}\right) +f_{k}\left(
f_{m+2}-f_{m+3}-f_{m+1}\right) =$\newline
$=2\left( f_{m}f_{k+1}-f_{k}f_{m+1}\right) =2\left( -1\right) ^{k}f_{m-k}.$%
\newline
The coefficient of $e_{3}$ is \newline
$f_{m+4}f_{k+6}-f_{m+3}f_{k+5}+f_{m+5}f_{k+2}-f_{m+2}f_{k+5}+f_{m+7}f_{k+1}-$
$f_{k+7}f_{m+1}=$\newline
$%
=f_{m}f_{k+2}-f_{m+2}f_{k}+f_{m+3}f_{k}-f_{m}f_{k+3}-f_{m+6}f_{k}+f_{m}f_{k+6}=
$\newline
$=f_{m}\left( f_{k+2}-f_{k+3}+f_{k+6}\right) +f_{k}\left(
-f_{m+2}+f_{m+3}-f_{m+6}\right) =$\newline
$=7\left( f_{m}f_{k+1}-f_{k}f_{m+1}\right) =7\left( -1\right) ^{k}f_{m-k}.$%
\newline
The coefficient of $e_{4}$ is\newline
$f_{m+5}f_{k+7}-f_{k+5}f_{m+7}+f_{m+6}f_{k+3}-f_{k+6}f_{m+3}+f_{m+1}f_{k+2}-$
$f_{m+2}f_{k+1}=$\newline
$%
=-f_{m}f_{k+2}+f_{k}f_{m+2}-f_{m+3}f_{k}+f_{k+3}f_{m}-f_{m}f_{k+1}+f_{k}f_{m+1}=
$\newline
$=f_{m}\left( -f_{k+2}+f_{k+3}-f_{k+1}\right) =0.$\newline
The coefficient of $e_{5}$ is\newline
$f_{m+6}f_{k+1}-f_{k+6}f_{m+1}+f_{m+7}f_{k+4}-f_{k+7}f_{m+4}+f_{m+2}f_{k+3}-$
$f_{k+2}f_{m+3}=$\newline
$%
=-f_{m+5}f_{k}+f_{k+5}f_{m}+f_{m+3}f_{k}-f_{k+3}f_{m}+f_{m}f_{k+1}-f_{k}f_{m+1}=
$\newline
$=f_{m}\left( f_{k+5}-f_{k+3}+f_{k+1}\right) +f_{k}\left(
-f_{m+5}+f_{m+3}-f_{m+1}\right) =$\newline
$=4\left( f_{m}f_{k+1}-f_{k}f_{m+1}\right) =4\left( -1\right) ^{k}f_{m-k}.$%
\newline
The coefficient of $e_{6}$ is\newline
$f_{m+7}f_{k+2}-f_{k+7}f_{m+2}+f_{m+1}f_{k+5}-f_{k+1}f_{m+5}+f_{m+3}f_{k+4}-$
$f_{k+3}f_{m+4}=$\newline
$=f_{m+5}f_{k}-f_{k+5}f_{m}-f_{m}f_{k+4}+f_{k}f_{m+4}-f_{m}f_{k+1}+$ $%
f_{k}f_{m+1}=$\newline
$=f_{m}\left( -f_{k+5}-f_{k+4}-f_{k-1}\right) +f_{k}\left(
f_{k+5}+f_{k+4}+f_{k-1}\right) =$\newline
$=-9\left( f_{m}f_{k+1}-f_{k}f_{m+1}\right) =-9\left( -1\right) ^{k}f_{m-k}.$%
\newline
The coefficient of $e_{7}$ is\newline
$f_{m+1}f_{k+3}-f_{k+1}f_{m+3}+f_{m+2}f_{k+6}-f_{k+2}f_{m+6}+f_{m+4}f_{k+5}-$
$f_{k+4}f_{m+5}=$\newline
$=f_{m}\left( -f_{k+2}+f_{k+4}+f_{k+1}\right) +f_{k}\left(
f_{m+2}-f_{m+4}-f_{m+1}\right) =$\newline
$=3\left( f_{m}f_{k+1}-f_{k}f_{m+1}\right) =3\left( -1\right) ^{k}f_{m-k}.$%
\newline
We obtain that\newline
$F_{k}\times F_{m}=\left( -1\right) ^{k}f_{m-k}\left(
-3e_{1}+2e_{2}+7e_{3}+4e_{5}-9e_{6}+3e_{7}\right) .$ \newline
Therefore\newline
$<F_{k}\times F_{m},F_{n}>=\left( -1\right) ^{k}f_{m-k}\left(
-3f_{n+1}+2f_{n+2}+7f_{n+3}+4f_{n+5}-9f_{n+6}+3f_{n+7}\right) =$\newline
$=-2f_{n+2}+2f_{n+1}+2f_{n}=0.$\medskip 

\begin{equation*}
\end{equation*}

\textbf{References}

\begin{equation*}
\end{equation*}%
\newline
\newline
[Akk; ] I Akkus, O Kecilioglu, \textit{Split Fibonacci and Lucas Octonions},
accepted in Adv. Appl. Clifford Algebras.\newline
[Br; ] R. Brown, A. Gray (1967), \textit{Vector cross products}, Commentarii
Mathematici Helvetici \qquad 42 (1)(1967), 222--236.\newline
[Fl, Sa; 15] C. Flaut, D. Savin, \textit{Quaternion Algebras and Generalized
Fibonacci-Lucas Quaternions}, accepted in Adv. Appl. Clifford Algebras 
\newline
[Fl, Sh; 13(1)] C. Flaut, V. Shpakivskyi,\textit{\ Real matrix
representations for the complex quaternions}, Adv. Appl. Clifford Algebras,
23(3)(2013), 657-671.\newline
[Fl, Sh; 13(2)] Cristina Flaut and Vitalii Shpakivskyi, \textit{On
Generalized Fibonacci Quaternions and Fibonacci-Narayana Quaternions}, Adv.
Appl. Clifford Algebras, 23(3)(2013), 673-688.\newline
[Fl, St; 09] C. Flaut, M. \c{S}tef\~{a}nescu, \textit{Some equations over
generalized quaternion and octonion division algebras}, Bull. Math. Soc.
Sci. Math. Roumanie, 52(100), no. 4 (2009), 427-439.\newline
[Gu;] I. A. Guren, S.K. Nurkan, \textit{A new approach to Fibonacci, Lucas
numbers and dual vectors}, accepted in Adv. Appl. Clifford Algebras.\newline
[Ha; ] S. Halici, On Fibonacci Quaternions, Adv. in Appl. Clifford Algebras,
22(2)(2012), 321-327.\newline
[Ha1; ] S Halici, \textit{On dual Fibonaci quaternions}, Selcuk J. Appl
Math, accepted.\newline
[Ho; 61] A. F. Horadam, \textit{A Generalized Fibonacci Sequence}, Amer.
Math. Monthly, 68(1961), 455-459.\newline
[Ho; 63] A. F. Horadam, \textit{Complex Fibonacci Numbers and Fibonacci
Quaternions}, Amer. Math. Monthly, 70(1963), 289-291.\newline
[Ja; 74] Nathan Jacobson (2009). \textit{Basic algebra I}, Freeman 1974 2nd
ed., 1974, p. 417--427.\newline
[Ke; ] O Kecilioglu, I Akkus, \textit{The Fibonacci Octonions,} accepted in
Adv. Appl. Clifford Algebras.\newline
[Ro; 96] M. Rost, \textit{On the dimension of a composition algebra}, Doc.
Math. J.,1(1996), 209-214.\newline
[Nu; ] S.K. Nurkan, I.A. Guren, \textit{Dual Fibonacci quaternions},
accepted in Adv. Appl. Clifford Algebras.\newline
[Sa, Fl, Ci; 09] D. Savin, C. Flaut, C. Ciobanu, \textit{Some properties of
the symbol algebras, Carpathian Journal of Mathematics}, 25(2)(2009), p.
239-245.\newline
[Si; 02] Z. K. Silagadze, \textit{Multi-dimensional vector product}, arxiv 
\newline
[Sc; 66] R. D. Schafer, \textit{An Introduction to Nonassociative Algebras},
Academic Press, New-York, 1966.\newline
[Sc; 54] R. D. Schafer, \textit{On the algebras formed by the Cayley-Dickson
process}, Amer. J. Math. 76 (1954), 435-446.\newline
[Sm; 04] W.D.Smith, \textit{Quaternions, octonions, and now, 16-ons, and
2n-ons; New kinds of numbers,\newline
}[Sw; 73] M. N. S. Swamy, \textit{On generalized Fibonacci Quaternions}, The
Fibonacci Quaterly, 11(5)(1973), 547-549.\newline
[Wo] http://mathworld.wolfram.com/FibonacciNumber.html

\begin{equation*}
\end{equation*}

Cristina FLAUT

Faculty of Mathematics and Computer Science,

Ovidius University,

Bd. Mamaia 124, 900527, CONSTANTA,

ROMANIA

http://cristinaflaut.wikispaces.com/

http://www.univ-ovidius.ro/math/

e-mail:

cflaut@univ-ovidius.ro

cristina flaut@yahoo.com

Vitalii SHPAKIVSKYI

Department of Complex Analysis and Potential Theory

Institute of Mathematics of the National Academy of Sciences of Ukraine,

3, Tereshchenkivs'ka st.

01601 Kiev-4

UKRAINE

http://www.imath.kiev.ua/

e-mail: shpakivskyi@mail.ru

\bigskip

\bigskip

\begin{equation*}
\end{equation*}

\end{document}